\documentclass[sigconf]{acmart}

\AtBeginDocument{%
  \providecommand\BibTeX{{%
    \normalfont B\kern-0.5em{\scshape i\kern-0.25em b}\kern-0.8em\TeX}}}

\copyrightyear{}
\acmYear{2021}
\acmDOI{}

\acmConference{}
\acmBooktitle{}
\acmPrice{}
\acmISBN{}

\usepackage{microtype}
\usepackage{tikz}
\usepackage{algorithm,comment,bm}
\usepackage[noend]{algpseudocode}
\usepackage{amsmath, graphicx, url, mdframed} 
\graphicspath{{figs/}}
\usepackage[normalem]{ulem}
\usepackage{pgf}
\newcommand{\R}{{\mathbf{R}}}

\newcommand{\cB}{\mathcal{B}} 

\newcommand{\cG}{\mathcal{G}}
 
\newcommand{\cN}{\mathcal{N}}

\newcommand{\subj}{\textnormal{subject~to}}

\newcommand{\iter}{k}

\newcommand{\xk}{x_{\iter}}
\newcommand{\wik}{w_{i,\iter}}
\newcommand{\rik}{r_{i,\iter}}
\newcommand{\xikp}{x_{i,\iter+1}}

\newcommand{\xik}{x_{i,\iter}}

\newcommand{\uikp}{u_{i,\iter+1}}
\newcommand{\uik}{u_{i,\iter}}

\newcommand{\yik}{y_{i,\iter}}

\newcommand{\yjk}{y_{j,\iter}}

\newcommand{\barx}{\bar{x}}

\newcommand{\baru}{\bar{u}}

\newcommand{\bary}{\bar{y}} 
 
\newcommand{\hG}{\hat{G}}

\newcommand{\T}{^\top}
\newcommand{\inv}{^{-1}}

\DeclareMathOperator*{\minimize}{minimize}

\newcommand\oprocendsymbol{\hbox{$\square$}}
\newcommand\oprocend{\relax\ifmmode\else\unskip\hfill\fi\oprocendsymbol}

\newtheorem{theorem}{Theorem}[section]

 \newtheorem{lemma}[theorem]{Lemma}
\newtheorem{remark}[theorem]{Remark} 
\newtheorem{example}[theorem]{Example}

\algnewcommand\init{\item[\textbf{Initialization:}]}
\algnewcommand\evol{\item[\textbf{Evolution:}]}

\copyrightyear{2022} 
\acmYear{2022} 
\setcopyright{acmlicensed}
\acmConference[AIES '22] {Proceedings of the 2022 AAAI/ACM Conference on AI, Ethics, and Society}{August 1--3, 2022}{Oxford, United Kingdom.}
\acmBooktitle{Proceedings of the 2022 AAAI/ACM Conference on AI, Ethics, and Society (AIES '22), August 1--3, 2022, Oxford, United Kingdom}
\acmPrice{15.00}
\acmDOI{10.1145/3514094.3534132}
\acmISBN{978-1-4503-9247-1/22/08}

\begin{document}
\fancyhead{}

\title{Achievement and Fragility of Long-term Equitability}

\author{Andrea Simonetto}
\email{andrea.simonetto@ensta-paris.fr}
\affiliation{%
  \institution{UMA, ENSTA Paris, Institut Polytechnique de Paris}
  \city{91120 Palaiseau}
  \country{France}
  }

\author{Ivano Notarnicola}
\email{ivano.notarnicola@unibo.it}
\affiliation{%
  \institution{Dept.~of Electrical, Electronic and Information Engineering, Università di Bologna}
  \city{40136 Bologna}
  \country{Italy}
  }

{\bf

\textcopyright {A. Simonetto, I. Notarnicola} {2022}. This is the author's version of the work. It is posted here for your personal use. Not for redistribution. The definitive Version of Record was published in {Proceedings of the 2022 AAAI/ACM Conference on AI, Ethics, and Society}, http://dx.doi.org/10.1145/3514094.3534132.
}
%
%
%
%
\begin{abstract}
Equipping current decision-making tools with notions of fairness, equitability, or other ethically motivated outcomes, is one of the top priorities in recent research efforts in machine learning, AI, and optimization. In this paper, we investigate how to allocate limited resources to {locally interacting} communities in a way to maximize a pertinent notion of equitability. In particular, we look at the dynamic setting where the allocation is repeated across multiple periods ({e.g., yearly}), the local communities evolve in the meantime ({driven by the provided allocation)}, and the allocations are modulated by feedback coming from the communities themselves. We employ recent mathematical tools stemming from data-driven feedback online optimization, by which communities can learn their {(possibly unknown)} evolution, satisfaction, as well as they can share information with the deciding bodies. We design dynamic policies that converge to an allocation that maximize equitability in the long term.  We further demonstrate our model and methodology with realistic examples of healthcare and education subsidies design in Sub-Saharian countries. One of the key empirical takeaways from our setting is that long-term equitability is fragile, in the sense that it can be easily lost when deciding bodies weigh in other factors (e.g., equality in allocation) in the allocation strategy. Moreover, a naive compromise, while not providing significant advantage to the communities, can promote inequality in social outcomes. 
\end{abstract}


\maketitle

\section{Introduction}

When limited resources have to be allocated to multiple parties, a natural question of fairness, equitability, or other ethically motivated notions, arises. To help tackling this problem, a number of tools from machine learning, AI, and optimization have appeared over the years. Different notions of outcome fairness have been proposed and, while they have been shown to be sometimes conflicting, they have proven their usefulness in fair decision-making. 

A more recent, and less studied, issue regards fair allocation in a dynamic setting, where the allocations are repeated over multiple time periods, and  they drive the evolution of the welfare of the parties. In this sense, initial allocations may have a lasting impact at future time periods and the evolution, i.e., how the funds are used to generate tangible welfare and social outcomes, has to be considered in designing allocations that are ethically motivated. 

A motivating example for this setting is the study of subsidies design in low-income countries to improve the population welfare and social indicators, like life expectancy, years of schooling, access to water, literacy, and many others. A multi-period funding campaign can drive these indicators up, but how the local communities use the subsidies is key in deciding how to allocate them to maximize pertinent notions of fairness. Other examples stem from grant allocation in academia, closing the gender gap by affirmative actions, designing subsidies for renewable and low-carbon energy sources, and {\sc covid} vaccine allocation.  

In this paper, we look at maximizing a novel notion of long-term equitability, which is constructed upon  connected neighboring communities’ outcomes. We consider the dynamic setting where the allocation is repeated across multiple periods ({e.g., yearly}), the local communities evolve in the meantime ({driven by the provided allocation)}, and the allocations are modulated by feedback coming from the communities themselves. We employ recent mathematical tools stemming from data-driven feedback online optimization, by which communities can learn their {(possibly unknown)} evolution, satisfaction, as well as they can share information with the deciding bodies. We design dynamic policies that converge to an allocation that maximize equitability in the long term.  We further demonstrate our model and methodology with realistic examples of subsidies design in Sub-Saharian countries.

The problem we are considering is challenging for a variety of reasons: \emph{(i)} the dynamics of the local communities are key in determining {how} the welfare is generated but it is generally unknown and can change over time: one needs a way to circumvent the need for learning the dynamical system and use only funding-to-social-outcome data to design fair allocations; which we provide. \emph{(ii)} The allocation must include community-driven feedback to trade off modeling errors and track time-varying dynamical conditions; which we include. \emph{(iii)} The policies must be easily explainable, they should democratically include local preferences, and compromise over different wanted outcomes should be easy to be made; which we provide. 

Among the significant takeaways of our empirical results is that long-term equitability is fragile, in the sense that it can be easily lost when deciding bodies weigh in other factors (e.g., equality in allocation) in the allocation strategy. Moreover, a naive compromise, while not providing significant advantage to the communities, can promote inequality in social outcomes. Finally, our results suggest that subsidies alone are not sufficient to drive equality, and they support the idea that investing in systemic changes is required.

\subsection{Related work}

Ethically motivated objectives have entered the mainstream in AI, machine learning, optimization, and decision-making, see for example the seminal works~\cite{Arrow1971, Rawls1971, Varian1973}, and the more recent \cite{Hardt2016, Calmon2017, Heidari2018, Heidari-Krause2019, Kallus2019, barocas-hardt-narayanan, Ben-Porat2021}. However, dynamic effects of the fair decisions on evolving populations is a far less studied area. Recently, a series of papers, among which~\cite{Hussein2019, Heidari2019, D'Amour2020, Creager2020, Heidari2021, Ge2021, Zhang2021, Wen2021, Chi2021}, have started to investigate these dynamic effects, e.g., by modeling decisions and dynamics as Markov Decision Processes and by using reinforcement learning as well as dynamic programming to design fair policies and algorithms. One of the main messages of these works is that fair decisions in a static context, may not be fair in a dynamic scenario, where populations and disadvantage groups evolve in response to the decisions taken at previous time periods. In addition, as carefully analyzed in~\cite{Hussein2019}, imposing some notion of fairness may drive unfairness in some other notion and, if the algorithms are naively designed, in the long-term the population may be globally worse off than when it started. 

Our paper finds similar conclusions of~\cite{Hussein2019}. It uses a theoretical-lighter approach than most of the aforementioned papers (online optimization instead of reinforcement learning), and it is able to infer the population dynamics by input-output data. In practice, our policies do not need to know the population dynamics, but only its long-term effects, together with the population's feedback on the imposed decisions. This makes our approach and policies easy to implement, interpret, and generalize. 

The idea to incorporate feedback in the decision-making process is not novel in general~\cite{Morik2020}, but here we use it in a specific way, that makes our policies converge to a long-term optimal solution. The technical tools we use stem from recent developments in online optimization with and without user's feedback, and in particular~\cite{Bianchin2021a,Bianchin2021,Coulson2019}. We also use the concept of preference elicitation and user's satisfaction from~\cite{Kahneman1979,Simonetto2021,Notarnicola2022}, which can be seen as a way to empower local communities and democratize the decision-making process, as also expressed in~\cite{Kasy2021}.

The metric we use is a novel notion of equitability inspired by~\cite{Hossain2020}, with the difference that is for us a cost violation based on neighborhood proximities. With this metric, we are able to formalize usual trade-offs between the will of the population, allocation equality, and equitability of social outcomes as points on a Pareto frontier. In particular, we are able to include the voice of the population in the decision-making process, as advocated in~\cite{Kasy2021,Yaghini2021}.  

Our running examples stem from subsidies design in low-income countries. Inspiring studies in this regard can be found in~\cite{imf2000, worldbank2005, LeBlanc2007, Oketch2016, Duflo2021} among others. However our method can be applied to other subsidies design as well, for instance health in richer countries~\cite{Chen2016}, energy subsidies~\cite{iea2021,Hortay2019,Line2019}, or {\sc{covid}} vaccine allocations~\cite{Matrajt2021, Liu2022}.

\section{A Dynamic Model}

\subsection{Funds-to-Welfare dynamics}

\begin{figure}
\centering
\includegraphics[width=0.4\textwidth]{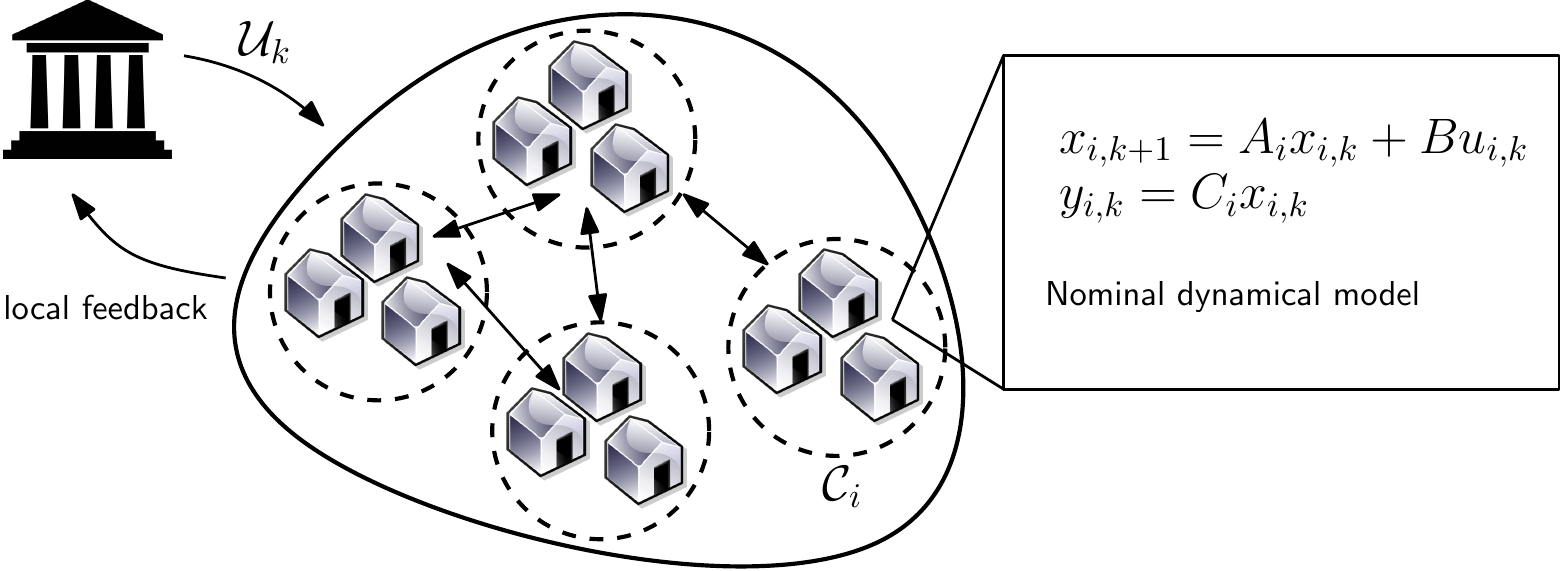}
\caption{A pictorial depiction of the problem setup consisting of a government agency that grants subsidies to local communities. The allocation is community-driven, based on maximizing equitability in the long term. All the notation is explained in the main text. }
\label{fig:setup}
\end{figure}

We consider a model for both a governmental funding agency $\mathcal{A}$ and the local communities $\mathcal{C}$, which we label $1, \ldots, N$. We also refer to Figure~\ref{fig:setup} for a pictorial representation of the main setup. 

Every funding period, say year $k$, the agency can allocate fundings $\mathcal{U}_k$ that are divided accordingly to some policy to the communities. We assume that each community can use the money they receive, say $\mathcal{U}_{i,k}$, to fund different activities (e.g., education, hospitals, infrastructure), and we label each activity as $1, \ldots, m$. In this way, the money that is allocated to community $i$ and activity $j$ in the $k$-th year, is indicated with $u_{ij,k}$. For simplicity, we stack all the $u_{ij,k}$ for the different activities in a vector $u_{i,k} \in \R^m_{+}$ (where $\R_{+}$ represents the non-negative real numbers).

Each community generates a welfare with the funding. A welfare can be the number of children who have received a scholarship, or the number of doctors who have been hired, and so forth. We capture the welfare of a given community $i$ at the year $k$ with the \emph{state} vector $x_{i,k} \in \R^n_{+}$, where $n$ are the number of dimensions that we consider (e.g., scholarships, doctors, ...). With funding $u_{i,k}$ for year $k$, a community generates welfare and this is captured by a dynamical equation: $x_{i,k+1} = f_i(x_{i,k}, u_{i,k})$, where $f_i:\R^n_{+}\times\R^m_{+} \to \R^n_{+}$ is a function that represents how the money is spent, and how the welfare is generated, and it is community specific.

One could specify different {functions $f$} depending on the knowledge of the local communities. One could even extend the present framework to probabilistic settings, where the dynamical equation is modeled via a Markov chain, either static or time-varying in time. For simplicity of exposition, we focus here on a linear relationship, such as $\xikp  = A_i\xik + B_i\uik + \wik$ for appropriate matrices $A_i\in \R^{n\times n}$, $B_i \in \R^{n\times m}$ and noise term $\wik {\in \R^{n}}$ taking care of possible errors in the model. We also explicitly avoid to consider local fundings in the equation, even though that could be incorporated into $\uik$ and properly factor out when asked for subsidies. 

\begin{example}\label{example}
Consider community $i$ receiving $20K$ USD for a new medical equipment ($u_{1,k}$), and $5K$ USD for buying children eLearning tools ($u_{2,k}$). We consider as welfare $x$ the number of people successfully diagnosed and cured $x_{1}$, and the number of children finishing that school year $x_{2}$. We could have, for example,
\begin{equation*}
\left[\begin{array}{c} x_{1,k+1} \\ x_{2,k+1} \end{array}\right] = \left[\begin{array}{cc} 0.8 & 0\\ 0 & 0.8 \end{array}\right]\left[\begin{array}{c} x_{1,k} \\ x_{2,k} \end{array}\right] + \left[\begin{array}{cc} 1.0 & 0\\ 0.5 & 1.0 \end{array}\right] \left[\begin{array}{c} u_{1,k} \\ u_{2,k} \end{array}\right],
\end{equation*}  
where matrix $A_i$ is modeled to take into account that medical and eLearning equipments degrade, so that if we do not invest in maintenance and new equipments, the number of people cured and the number of children attending school is going to decrease. We also model in $B_i$ the interaction between the number of people cured (some of which may be children), and the number of children finishing the school year. This last interaction between investment in healthcare having effect into schooling is well known~\cite{dataset-health-1,dataset-edu-1}. 
\end{example}

In the following, we assume that the linear dynamical system is \emph{stable}, meaning that the eigenvalues of $A_i$ are all inside the unit circle. This is a well-motivated assumption, since typically communities cannot generate more welfare without (local and external) fundings, and with no fundings, the welfare will eventually decrease to the stable zero equilibrium.   

Once welfare is generated, one can measure the satisfaction of a community based on different (standard) indicators. For example, one could measure or estimate the life expectancy at birth, and/or the number of years of schooling, or how many miles one has to travel to get access to drinkable water. Another standard indicator is the Human development index~\cite{hdi}. We leave the freedom here to define a vector of indicators $y_{i,k} \in \R^{p}_{+}$, which we assume proportional to the welfare via a suitable matrix, i.e., $y_{i,k} = C_{i} x_{i,k} + \rik$, with the addition of a noise term $\rik$ that captures errors in the linear models, as well as errors in measuring the indicator\footnote{More complex equations where $\uik$ also plays a role in the satisfaction can also be considered.}. Once again, the linear relationship is a reasonable assumption based on our current understanding~\cite{dataset-health-1,dataset-edu-1}. 

Putting our model together, we will consider a funding-to-welfare dynamical model with a measure of satisfaction as
\begin{equation}\label{noisemodel}
\mathcal{C}_i : \left\{\begin{array}{rcl}
  \xikp & = &A_i\xk + B_i\uik + \wik
  \\
  \yik & = & C_i\xik + \rik
  \end{array}\right.
\end{equation}

The matrices $(A_i, B_i, C_i)$ are generally unknown to the funding agency, but can be estimated locally via data-driven approaches from historical data or a priori knowledge, as we explain in Section~\ref{sec:learning}. Once this triplet is known, under the stability assumption, one can derive the notion of nominal equilibrium, meaning which funding $u_{i,k}$ will maintain the current level of welfare and satisfaction in the long term, without noise. This nominal equilibrium triplet $(\barx_i,\baru_i,\bary_i)$ is found by solving the linear system:
\begin{equation}\label{eq:equilibrium}
\mathcal{C}_i(\barx_i,\baru_i,\bary_i) : \left\{\begin{array}{rcl}
  \barx_i & = & A_i\barx_i + B_i\baru_i
  \\
  \bary_i & = & C_i\barx_i
  \end{array}\right.
\end{equation}
and in particular, under the stability assumption, $\bary_i = C_i  (I_n-A_i)\inv B_i \baru_i$. We let $G_i = C_i (I_n-A_i)\inv B_i$ for the following, and we call it the input-output static map. 

Our main aim is to design a policy to allocate $\uik$ from year to year, such that eventually we maximize the long-term satisfaction $\bary_i$, with a pertinent notion of equitability. We see how we define the latter next. 


\subsection{Optimizing for the long-term}

To determine a suitable policy sequence of allocations $\uik$, one has to define a metric to optimize over. First, we define a graph $\cG$ connecting close-by communities. Our intention is to derive a metric of equitability that depends on local communities and on how these communities see themselves \emph{with respect to} other nearby communities. 

Let then $\cG$ be a graph, such that the nodes are the communities and the edges are the links between them. In this way, if community $i$ is close-by community $j$ (i.e., $j$ is one of its neighbors), there is a link between them. We call $\cN_i$ the set of neighbors of community $i.$ We further call {$N_i$} the number of neighbors community $i$ has.

Second, we define equitability. In standard literature, equitability would be defined as imposing that $y_{i,k} = y_{j,k}$ for each pair of local communities $i,j$. This is a too stringent requirement in many applications, and it has already been relaxed into groups, and/or equitability violations~\cite{Hossain2020}. We follow the route of relaxing this constraint first by asking that $y_{i,k} = \frac{1}{N_i} \sum_{j \in \cN_i}  y_{j,k}$, for all $i$'s (meaning looking at group averages within the neighborhood, which is what most communities have access to), and then moving the constraint in the cost function so that its violation is properly penalized. In particular, we consider a neighborhood equitability-violation metric (NEqM) as follows,
\begin{equation}\label{NEqM}
\textrm{(NEqM)} \qquad \psi_i(\yik, \{\yjk\}_{j \in \cN_i}) = \Big\| \, \yik - \frac{1}{N_i}  \sum_{j \in \cN_i}\yjk \,\Big\|^2.  
\end{equation}
Convex function $\psi_i(\yik, \{\yjk\}_{j \in \cN_i})$ represents the distance of the local community $i$ from perfect equitability with its neighboring communities. 

Finally, we can define our optimization strategy which is captured by the optimization problem:
\begin{align}\label{eq:centralized_pb}
\begin{split}
  \minimize_{\{\baru_i,\barx_i, \bary_i\}_{i \in [1,N]}} \: & \: \varphi(\{\baru_i\}_{i \in [1,N]}) + \sum_{i=1}^N \psi_i(\bary_i, \{\bary_j\}_{j \in \mathcal{N}_i}) 
  \\
  \subj \: & \: \left\{\begin{array}{c}\textrm{Eq. } \eqref{eq:equilibrium} \: \textrm{for all communities } i \\
  \baru_i \geq 0, \quad \sum_{i=1}^N \baru_i \leq s_{\max},\end{array}\right.
\end{split}
\end{align}
where, for completeness, we have also added a global convex cost $\varphi(\{\baru_i\}_{i \in [1,N]})$ that could impose government preferences on the division of allocations among the different activities and among different communities. Note also the global budget constraint: $\sum_{i=1}^N \baru_i \leq s_{\max}$ for a maximal funding $s_{\max}$, that is imposed not to run over-budget. 

By exploiting the map $G_i$, Problem~\eqref{eq:centralized_pb} can be equivalently recast into the 
following \emph{funding-only} optimization problem
\begin{align}\label{eq:centralized_pb_reduced}
\begin{split}
  \minimize_{\{\baru_i\}_{i \in [1,N]}}  \: & \: \varphi(\{\baru_i\}_{i \in [1,N]}) + \sum_{i=1}^N \psi_i(G_i\baru_i, \{G_j\baru_j\}_{j \in \cN_i}),\\ \subj \:  &  \: \baru_i \geq 0, \quad \sum_{i=1}^N \baru_i \leq s_{\max}.
  \end{split}
\end{align}

In principle, if the maps $G_i$ for all $i$ were known accurately, the government could allocate subsidies in a long-term equitable way, and its policy would be the same at every funding rounds. We call this strategy, the static open loop (SOL) policy. In reality, the dynamics of the local communities, and therefore $G_i$, are known only approximately, the dynamics may change in time, and therefore one has to incorporate the \emph{feedback} of the local communities, while solving~\eqref{eq:centralized_pb_reduced}. 


\section{Policies}
\subsection{Static ideal policy}

Incorporating feedback in the decision-making process is key in delivering high societal-value outcomes, in the face of uncertainties. We use here recent tools stemming from data-driven and feedback-based online optimization. 

First, however, we look at a static policy. With the aim of solving~\eqref{eq:centralized_pb_reduced}, we set up a gradient iteration, as follows in the SOL policy. 
\medskip

\begin{mdframed}
\textbf{SOL Policy}
\begin{enumerate}
\item Start with a tentative $u_{i}[0]$ for all communities $i$, a choice of stepsize $\gamma>0$ and a choice of maximum iteration steps $\ell_{\max}$.

\item For each funding period $k$:

\begin{itemize}
\item For all $\ell \in[0, \ell_{\max}]$ and all communities $i$, iterate with a gradient step:
\begin{equation*}\qquad\quad\left\lfloor
\begin{array}{l}
 \hat{u}_i[\ell] = u_i[\ell]
    - \gamma \Big(
    \nabla_{u_i}\varphi(\{u_i[\ell]\}_{i \in [1,N]}) + \\ 
   \qquad\qquad \sum_{i=1}^N G_i^\top \nabla_{u_i}\psi_i(G_iu_i[\ell], \{ {G_j} u_j[\ell]\}_{j \in \cN_i})
    \Big), \\
   \{u_i[\ell+1]\}_{i \in [1,N]} = \Pi_{\cB}(\{\hat{u}_i[\ell]\}_{i \in [1,N]})
\end{array}\right.
\end{equation*}
where $\Pi_{\cB}(\cdot)$ is the projection onto the convex set $\cB = \left\{\{u_i\}_{i \in [1,N]} \mid u_i \geq 0, \sum_{i=1}^N u_{i} \leq s_{\max}\right\}$. 
\item Set and implement $u_{i,k} = u_i[\ell_{\max}+1]$
\end{itemize}
\end{enumerate}
\end{mdframed}
\medskip

We know, by standard results in convex analysis and projected gradient method~\cite{Nesterov2004, Nocedal2006}, that the SOL policy delivers a sequence $\{u_i[\ell]\}$ that converges to the solution of Problem~\eqref{eq:centralized_pb_reduced} for sufficiently small stepsize $\gamma$, as captured in the following lemma. 

\begin{lemma}\label{lemma:1}
Consider problem Problem~\eqref{eq:centralized_pb_reduced} and its optimizers $\baru_i^*$. Assume function $\varphi$ to be convex. Assume also that the cost function $\varphi(\{\baru_i\}_{i \in [1,N]}) + \sum_{i=1}^N \psi_i(G_i\baru_i, \{G_j\baru_j\}_{j \in \mathcal{N}_i})$ is $L$-smooth (i.e., it has a $L$-Lipschitz continuous gradient). Then, choosing $\gamma < 2/L$, the SOL policy will deliver a sequence for which $\lim_{\ell \to \infty} \|u_i[\ell] - \baru_i^*\| = 0$ for all communities $i$.
\end{lemma}

Lemma~\ref{lemma:1} ensures that if $\ell_{\max}$ is taken sufficiently large, we can set $\baru_i = u_{i,k}$. This allocation will be the same across the years, since $G_i$'s are not changing in the SOL policy model. This is therefore the best static long-term policy\footnote{Under the smoothness requirements of the lemma, convergence can also be accelerated via Nesterov's algorithm, but we do not explore this here.}.   

\subsection{Going beyond the utopia}

It is not hard to see that the presented SOL policy is a utopian goal. Not only the triplet $(A_i, B_i, C_i)$ is not known (and therefore $G_i$ is unknown), but this triplet is only an idealized model of an underlying more complex system. In this context, with the SOL policy we would like reality to converge to an ideal model that does not exist in practice. In fact, it is much more reasonable to aim at the goal of solving the problem, 
\begin{align}\label{eq:centralized_pb_reduced_reason}
\begin{split}
  \minimize_{\{\baru_i\}_{i \in [1,N]}}  \: & \: \mathbb{E}\Big[\varphi(\{\baru_i\}_{i \in [1,N]}) + \\ & \!\!\!\sum_{i=1}^N \psi_i(G_i\baru_i + H_i \wik + \rik, \{G_j\baru_j+ H_j w_{j,k} + r_{j,k}\}_{j \in \cN_i})\Big], \\ \subj \: & \: \bar{u}_i \geq 0, \quad \sum_{i=1}^N \baru_i \leq s_{\max},
  \end{split}
\end{align}
where we have set $H_i = C_i(I_n-A_i)^{-1}$, the expectation $\mathbb{E}[\cdot]$ is with respect to all the random variables, and we have reintroduced the modeling errors $\wik$, as well as the measurement errors $\rik$ (Cf. Eq.~\eqref{noisemodel}). This is not a static optimization problem, since the random variable distribution may change in time, so the best policy will change at every funding period.

Even when looking at the modified~\eqref{eq:centralized_pb_reduced_reason}, several challenges are still present. First, the input-output static maps $G_i$ are only known up to a certain accuracy, say $\hG_i$, due to the noise terms $\wik, \rik$, and they may drift over the years (since $A_i, B_i, C_i$ may do that). This requires a modification of the gradient descent to incorporate actual community feedback. One could take for example feedback $y_i[\ell]$ for allocation $u_i[\ell]$; we know that at steady state $y_i[\ell] = G_i u_i[\ell]+ H_i \wik + \rik$, and therefore we could take a modified descent as:
\begin{multline}\label{theproblem}
{\hat{u}_i[\ell]} = u_i[\ell]
    - \gamma \Big(
    \nabla_{u_i}\varphi(\{u_i[\ell]\}_{i \in [1,N]}) + \\ \sum_{i=1}^N \hG_i\T \nabla_{u_i}\psi_i(y_i[\ell], \{y_j[\ell]\}_{j \in \mathcal{N}_i})
    \Big), \quad \forall i.
\end{multline}
In particular, we take the feedback $y_i[\ell]$ from the community due to the input $u_i[\ell]$. Two important remarks are now necessary. First, we cannot expect to be asking communities feedback very frequently, for many logistic reasons and also because we need time before the funding will create a change in welfare. So we will let $\ell = k$, in a way that we run one gradient per funding period. Second, since one can expect that the gradient may be slow at converging, we warm start it by setting the initial condition equal to the solution of the SOL policy.

With this in place, we can devise our first dynamic close loop (DCL) policy as follows.

\medskip

\begin{mdframed}
\textbf{DCL Policy}
\begin{enumerate}
\item Start with $u_{i,0}$ for all communities $i$, being equal to the SOL policy solution with an estimated $\hG_i$, a choice of stepsize $\gamma>0$. 

\item For all funding period $k$ and all communities $i$, iterate with a gradient step:
\begin{itemize}
\item Ask for noisy community feedback $\yik = C_i \xik + \rik$.

\item Calculate:
\begin{equation*}\quad\left\lfloor
\begin{array}{l}
 \hat{u}_{i,k} = u_{i,k}
    - \gamma \Big(
    \nabla_{u_i}\varphi(\{u_{i,k}\}_{i \in [1,N]}) + \\
    \qquad\qquad \sum_{i=1}^N \hG_i\T \nabla_{u_i}\psi_i(\yik, \{\yjk\}_{j \in \cN_i})
    \Big), \\
   \{\uikp\}_{i \in [1,N]} = \Pi_{\cB}(\{\hat{u}_{i,k}\}_{i \in [1,N]})
\end{array}\right.
\end{equation*}
where $\Pi_{\cB}(\cdot)$ is the projection onto the convex set $\cB = \left\{\{u_i\}_{i \in [1,N]} \mid u_i \geq0, \sum_{i=1}^N u_{i} \leq s_{\max}\right\}$. 

\item Implement $\uikp$. 
\end{itemize}
\end{enumerate}
\end{mdframed}
\medskip

Notice immediately that the feedback is given on the evolving dynamical system $\yik = C_i \xik + \rik$, as we do not have a steady-state equilibrium point in general (so $y_{i,k} \neq G_i u_{i,k}+ H_i \wik + \rik$). However this is far from detrimental, since it allows for correcting for modeling errors. In addition, under reasonable assumptions, we know that for small enough stepsizes $\gamma$, the DCL policy will deliver a sequence $u_{i,k}$ that converges to the solution of Problem~\eqref{eq:centralized_pb_reduced_reason} within an arbitrarily small error bound, as follows.

\begin{lemma}\label{lemma:2}
Consider problem Problem~\eqref{eq:centralized_pb_reduced_reason} and its time-varying optimizers $\baru_i^{*,k}$. Assume function $\varphi$ to be strongly convex. Assume also that the cost function $\varphi(\{\baru_i\}_{i \in [1,N]}) + \sum_{i=1}^N \psi_i(G_i\baru_i + H_i \wik + \rik, \{G_j\baru_j+ H_j w_{j,k} + r_{j,k}\}_{j \in \cN_i})$ is $L$-smooth (i.e., it has a $L$-Lipschitz continuous gradient). Then, there exists a sufficiently small $\gamma$ for which the DCL policy will deliver a sequence for which $\limsup_{k \to \infty} \mathbb{E}[\|u_{i,k} - \baru_i^{*,k}\|] = E < \infty$ for all communities $i$.

In particular, the error bound $E$ is directly proportional to how the time-varying optimizers $\baru_i^{*,k}$ change in time, the approximation error $\|G_i - \hG_i\|$, and the single-point gradient approximation error for the expectation. 
\end{lemma}

Lemma~\ref{lemma:2} is an embodiment of Theorem~{6.1} and Proposition~4.2 of~\cite{Bianchin2021}, where the various proportionality constants in the error are spelled out. The interesting thing here is that the DCL policy can successfully take into account community-driven information to shape the decision-making process and deliver near-optimal allocation. The asymptotical error is proportional to how accurately we know and can model the underlying dynamics. 

Before moving on, it is interesting to take the time to analyze the DCL policy iteration once more.  
Note that $u_{i,k}$ is initialized with the SOL policy solution: it is the ideal long-term strategy that government could aim at. This ideal scenario is then faced with short-term reality in terms of feedback terms $\yik$, which is the voice of the community. All is then arranged together and projected into the allowed budget $s_{\max}$. {\bf The convergence lemma is then the proof that combining an idealized model with pertinent feedback can work in delivering long-term optimal allocations.}  

\subsection{Learning and re-learning}

Since Problem~\eqref{eq:centralized_pb_reduced_reason} has already time-varying optimizers $\baru_i^{*,k}$, one can imagine to modify the DCL policy, introducing a re-learning of $\hG_i$, whenever it is required, or whenever new data becomes available. As we express below. 

\medskip

\begin{mdframed}
\textbf{DCL Policy with re-learn (DCL+)}
\begin{enumerate}
\item Start with $u_{i,0}$ for all communities $i$, being equal to the SOL policy solution with an estimated $\hG_i$, a choice of stepsize $\gamma>0$. 

\item For all funding period $k$ and all communities $i$, iterate with a gradient step:
\begin{itemize}
\item Ask for noisy community feedback $\yik = C_i \xik + \rik$.

\item Calculate:
\begin{equation*}\quad\left\lfloor
\begin{array}{l}
 \hat{u}_{i,k} = u_{i,k}
    - \gamma \Big(
    \nabla_{u_i}\varphi(\{u_{i,k}\}_{i \in [1,N]}) + \\
    \qquad\qquad \sum_{i=1}^N \hG_i\T \nabla_{u_i}\psi_i(\yik, \{\yjk\}_{j \in \cN_i})
    \Big), \\
   \{\uikp\}_{i \in [1,N]} = \Pi_{\cB}(\{\hat{u}_{i,k}\}_{i \in [1,N]})
\end{array}\right.
\end{equation*}
where $\Pi_{\cB}(\cdot)$ is the projection onto the convex set $\cB = \left\{\{u_i\}_{i \in [1,N]} \mid u_i \geq0, \sum_{i=1}^N u_{i} \leq s_{\max}\right\}$.

\item Implement $\uikp$. 
\item Add the new data $(u_{i,k}, x_{i,k}, y_{i,k})$ to the historical data and re-learn $\hG_i$.
\end{itemize}
\end{enumerate}
\end{mdframed}
\medskip

This new policy will converge very similarly to the DCL policy, since the sources of errors are the same, but it has the advantage to incorporate new information concurrently as the implementation of new allocations. 

\subsection{Learning $G_i$}
\label{sec:learning}

The above described policies rely on leaning of the map $G_i$ for each community, which can be obtained from historical data. We describe here briefly three ways to learn $G_i$. 

The most direct (somewhat naive way) is to learn $G_i$ with $y_i$ and $u_i$ data with linear regression, i.e., by fitting a line on the data points. This discards the underlying dynamical system, but can be an effective way, especially if the dynamics is fast (i.e., when all the eigenvalues of $A_i$ are close to $0$). 

A more appropriate way to incorporate dynamics in the learning is via a behavioral approach, which involves the solution of a system {of linear equations}. The {mathematical} details can be found in~\cite{Bianchin2021a,Bianchin2021,Coulson2019}, and they require {sufficiently informative data to work.}

The most sophisticated way to learn $G_i$, is to perform system identification on $u_i, x_i$ and reconstruct all the matrices $A_i, B_i, C_i$ and then set $G_i = C_i (I_n - A_i)^{-1} B_i$. While this approach could estimate the whole dynamical system, it is often an overkill (since we only need $G_i$, and $x_i$ maybe not very easy to estimate). For the interested reader, we refer to~\cite{Ljung1999}. 

\section{Cost Choices}

\subsection{Other fairness metrics: worst-case}

While we have studied in more details an equitability violation metric as defined in~\eqref{NEqM}, one can substitute it with something more pertinent to specific situations at hand. For example, a concept that has received attention is to minimize in the worst-case scenario, e.g.~\cite{Emily2021}, which can be interpreted as Rawlsian ``maximin'' fairness~\cite{Rawls1971}. In this case, one could consider the convex metric:
\begin{equation}\label{WC-NEqM}
  \textrm{(WC-NEqM)} 
  \qquad \psi_i(\yik, \{\yjk\}_{j \in \cN_i}) 
  = \Big\| \, \yik - \frac{1}{N_i}  \sum_{j \in \cN_i}\yjk \,\Big\|_{\infty}.  
\end{equation}
With {this different cost}, most of our theoretical discussions still hold, with the exception that {one would need to adapt the optimization method to handle non-differentiable functions.}

\subsection{Design of the funding cost}\label{sec:design}

So far we have focused on the design of the output-dependent cost $\sum_{i=1}^N \psi_i(\bary_i, \{\bary_j\}_{j \in \mathcal{N}_i})$ in Problem~\eqref{eq:centralized_pb_reduced_reason}, leaving to the deciding bodies the construction of a reasonably motivated funding-dependent cost $\varphi(\{\baru_i\}_{i \in [1,N]})$. We will see in the examples how this choice is not trivial and can jeopardize in practice the equitability that we have tried to enforce with the output-dependent cost. But before that, we examine here potential cost design. 

The baseline design for $\varphi(\{\baru_i\}_{i \in [1,N]})$ is the weighted equal allocation\footnote{This can be referred to as allocation parity, equal partitioning, egalitarian allocation, and so forth.}:
\begin{equation}
\varphi(\{\baru_i\}_{i \in [1,N]}) =  \sum_{i,j \in [1,N] }\Big\| \, \baru_i - \baru_j \,\Big\|^2,
\end{equation}
enforcing the allocation of the same funding to every community, possibly weighted by population. This is a commonly advocated allocation, often believe ethically motivated and ``fair''. This is, among many examples, the strategy the European Union has adopted to allocate {\sc{covid}} vaccines over multiple time periods to its member's countries.

This allocation and cost is at odds with our equitability metric. In fact, since at equilibrium $\bar{y}_{i} = G_i \bar{u}_i$, and for the equal allocation $\bar{u}_i = \bar{u}_j$ for all $i, j$, then the equitability metric becomes determined by the differences in $G_i$ among the different communities and cannot be decreased. For example $\bar{y}_{i} - \bar{y}_{j} = (G_i - G_j)\bar{u}_i$. In fact, as we will see, if the resources augment ($s_{\max}$ augments), then $\bar{u}_i$ augments and therefore the equitability decreases. In this context, allocating equally among different communities (even factoring in population size) can foster inequalities in the outcome.

\subsection{Power and democracy}

As discussed, e.g., in~\cite{Kasy2021}, the choice of objective functions is intimately connected with the political economy question of who has ownership and control rights over data and algorithms, and how they drive the welfare of the people who have not chosen it. In fact, the decisional ``power'' resides in whom designs the objective function. Inspired by~\cite{Simonetto2021,Notarnicola2022}, one could enlarge the decision-making process to make it as inclusive as possible by considering an additional satisfaction term in the objective. The idea is to incorporate the satisfaction of a particular decision learned based on the community feedback on it. 

Let us imagine that each community has a function that measures their dissatisfaction for a particular decision and/or outcome \emph{relative} to their neighbors\footnote{Other absolute functions can also be considered in general.}. Let such functions be $\Delta_i$ defined as a monotone function of input and output discrepancies:
\begin{equation}
\Delta_i = \Delta_i\Big(\Big\| \, \bar{u}_i - \frac{1}{N_i}  \sum_{j \in \mathcal{N}_i} \bar{u}_j \,\Big\|^2, \Big\| \, \yik - \frac{1}{N_i}  \sum_{j \in \mathcal{N}_i}\yjk \,\Big\|^2\Big)
\end{equation}
Functions $\Delta_i$ are unknown but can be learned as explained in~\cite{Simonetto2021,Notarnicola2022} by leveraging surveys or other types of feedback. 

Putting this together, one could consider the composite cost function for Problem~\eqref{eq:centralized_pb} as  
\begin{multline}\label{demos}
\bar{f} := \overbrace{\underbrace{ \!\!\!\sum_{i,j \in [1,N] }\!\!\varrho\Big \| \, \bar{u}_i - \bar{u}_j \,\Big\|^2 + \sum_{i \in [1,N] }\Big\| \, \bar{y}_i - \frac{1}{N_i}  \sum_{j \in \mathcal{N}_i}\bar{y}_j \,\Big\|^2}_{\textrm{what the funding agency wants}}}^{(A) \quad\qquad \textrm{vs.} \qquad (B)}+ \\ \underbrace{ \!\!\sum_{i\in [1,N]}\sigma \Delta_i\Big(\Big\| \, \bar{u}_i - \frac{1}{N_i}  \sum_{j \in \mathcal{N}_i} \bar{u}_j \,\Big\|^2, \Big\| \, \bar{y}_i - \frac{1}{N_i}  \sum_{j \in \mathcal{N}_i}\bar{y}_j \,\Big\|^2\Big)}_{\textrm{what the local communities want}}
\end{multline}
with $\rho, \sigma \geq 0$ weighting terms. 
Equation~\eqref{demos} represents the whole democratic decision process, trading-off funding agencies' wishes with people's wishes, and on a finer scale, allocation equality $(A)$ with social outcome equality $(B)$.  

We are now ready to explore all we have discussed so far in two numerical examples.

\section{Examples}

\subsection{Health funding allocations in Sub-Saharan countries}

\subsubsection{Setting}

Our first example portraits the allocation of health fundings to selected countries in Sub-Saharan Africa. We put ourselves in a realistic (yet fictitious) setting of a global humanitarian organization who is deciding how to divide a maximum budget among several countries, in order to increase people's life expectancy at birth. We make a few simplifying assumptions for the sake of clarity: we assume that health funding is responsible for life expectancy~\footnote{While correlation has been empirically observed, causality has not been proven in general, but likely in low-income countries, like the ones considered here, see~\cite{dataset-health-1}.}, and that all the health funding (on top of the funding at time zero) comes from the humanitarian organization without local contribution. These two assumptions are made to keep the models simple and yet to capture the main issues. One could remove the assumptions by making the models more complex, but we do not pursue this here. 

\begin{figure}[b]
\centering
\includegraphics[width = 0.45\textwidth]{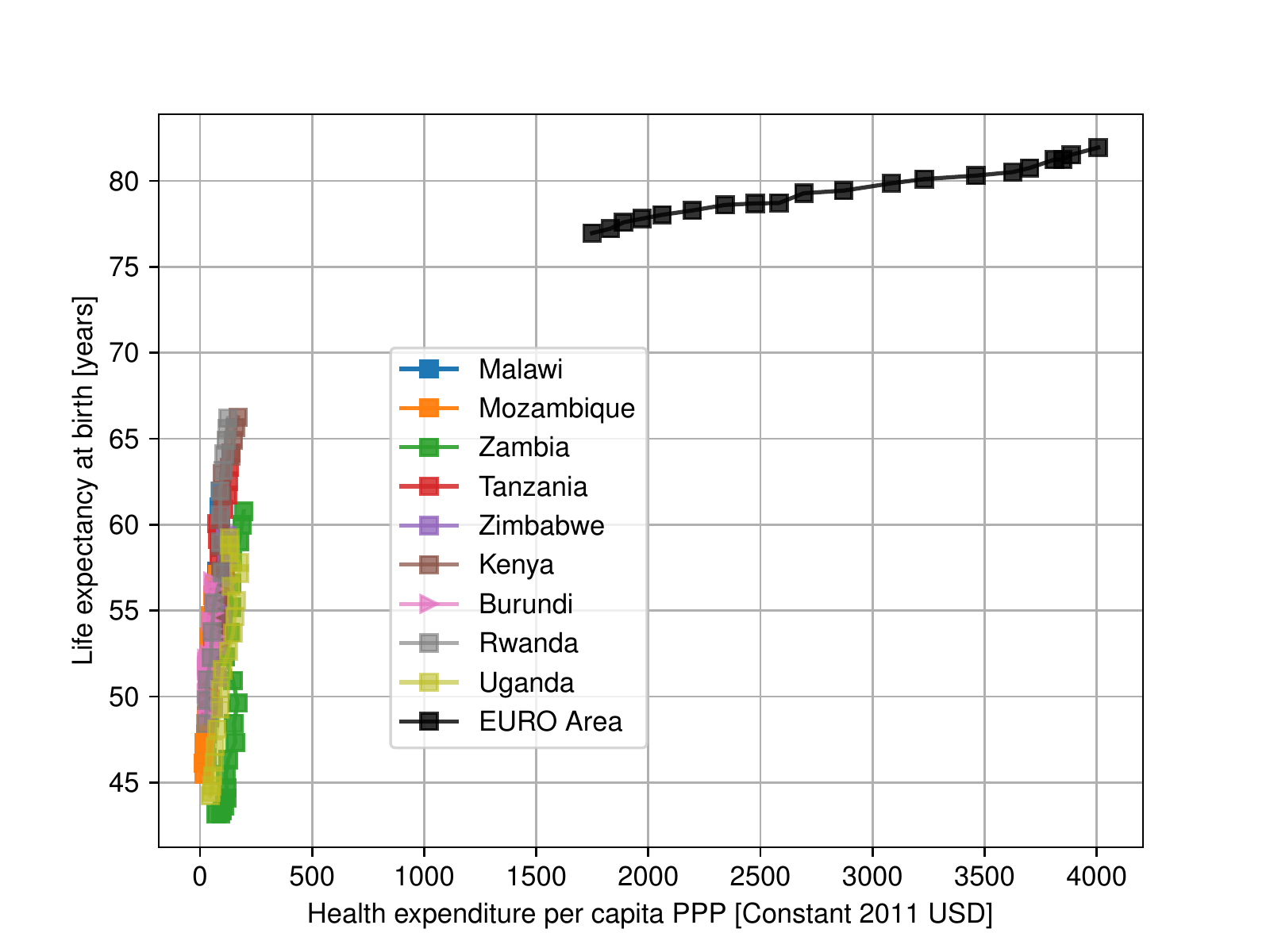}
\vskip-4.25cm\hspace*{3.5cm}\includegraphics[width = 0.2\textwidth]{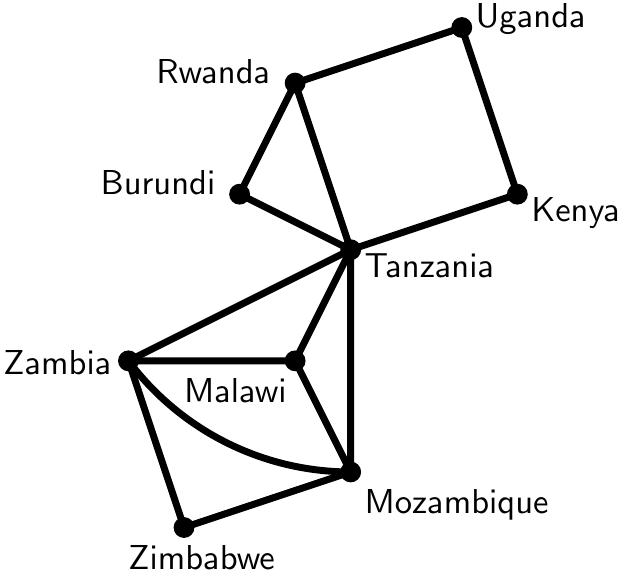}
\vskip1cm
\caption{The data available for the first example from~\cite{dataset-health-2}, with a pictorial depiction of how the countries are connected. The data represents the health expenditure vs.~the life expectancy at birth from 1985 to 2015. }
\label{fig.countries}
\end{figure}

We consider nine ({$N=9$}) countries, whose neighborhoods can be inferred by geography and pictorial represented in Figure~\ref{fig.countries}. We also collect health expenditure per capita PPP\footnote{Purchasing Power Parity.} (at constant 2011 USD) and the expected life expectancy at birth from~\cite{dataset-health-1, dataset-health-2}, from 1985 till 2015, see Figure~\ref{fig.countries}. Compared to the {\sc{euro}} Area, we see how the selected countries have a similar behavior. 

We estimate $G_i$ for all countries with linear regression from the data and we compute sensible $A_i, B_i, C_i$ to match the seen data. In particular, we set $B_i = 1$, while $A_i = 0.5 + \epsilon_i$ ($\epsilon_i$ being a zero-mean Gaussian noise term with standard deviation $0.1$) and $C_i = G_i (I-A_i)$. Having more or different data, one could also have done the estimation in a different fashion. 

\begin{figure*}
\centering
\includegraphics[width = 0.75\textwidth]{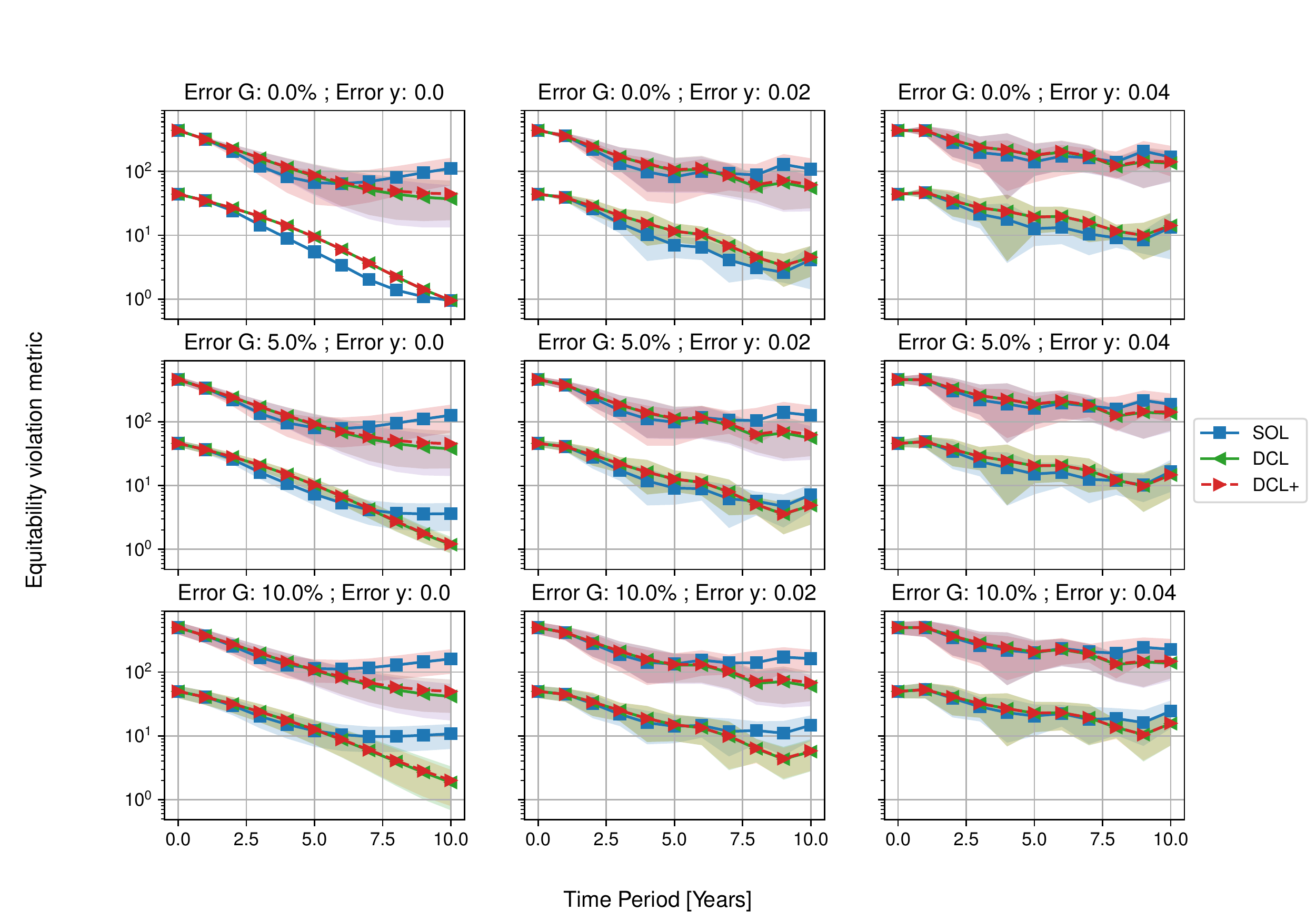}
\caption{Equitability violation metric for the different scenarios considered in the health subsidies example. These figures show that if the dynamical model is not accurately estimated, then both DCL and DCL+ policies, which integrate community feedback, are performing the best.}
\label{fig-2}
\end{figure*}

\subsubsection{Scenarios}

With this in place, we look at funding allocation for a ten-year time period. We set $s_{\max}$ to be a function of time that increases linearly every year to reach a $50\%$ increase of the nominal budget after $10$ years (this is to mimic the natural occurrence in the past 20 years). We also impose the constraints that the {entire} budget has to be used, and no country can get less than the initial health expenditure. Finally, we consider the cost function,
\begin{equation}\label{cost-ex1}
f_{\textrm{ex.1}} =  \!\!\sum_{i,j \in [1,N] }\varrho\Big \| \, \bar{u}_i - \bar{u}_j \,\Big\|^2 + \sum_{i\in [1,N] } \Big\| \, \bar{y}_i - \frac{1}{N_i}  \sum_{j \in \mathcal{N}_i}\bar{y}_j \,\Big\|^2
\end{equation}
with $\varrho = 0$, to focus on long-term equitability alone (see Eq.~\eqref{NEqM}).  

We consider then different scenarios, 
\begin{itemize}
\item[(S1)] The nominal case, where we have no errors and the dynamic is exactly reconstructed. 
\item[(S2)] A series of noisy cases, where the feedback on $y_{i,k}$ is given with some multiplicative\footnote{Meaning $y_{i,k} = (1+r_{i,k})\times(\textrm{true value})$. Multiplicative noise slightly changes our model, without affecting the theory, and it is more adapted here.} Gaussian noise $r_{i,k}$, with zero-mean and standard deviation: $\{0.02, 0.04\}$. 
\item[(S3)] A series of erroneously estimated $G_i$, whereby $G_i$ is estimated wrongly with a $\{5\%, 10\%\}$ standard deviation error.  
\item[(S4)] A series of evolving dynamics scenarios, where $G_i$ is estimated with some error level at the beginning, but then it drifts to other values (so the estimation gets less and less accurate as time progress). We model this as a slow drift of $G_i$ towards the $G_i$ of the {\sc euro} Area.  
\end{itemize}

\begin{figure*}
\centering
\includegraphics[width = 0.75\textwidth]{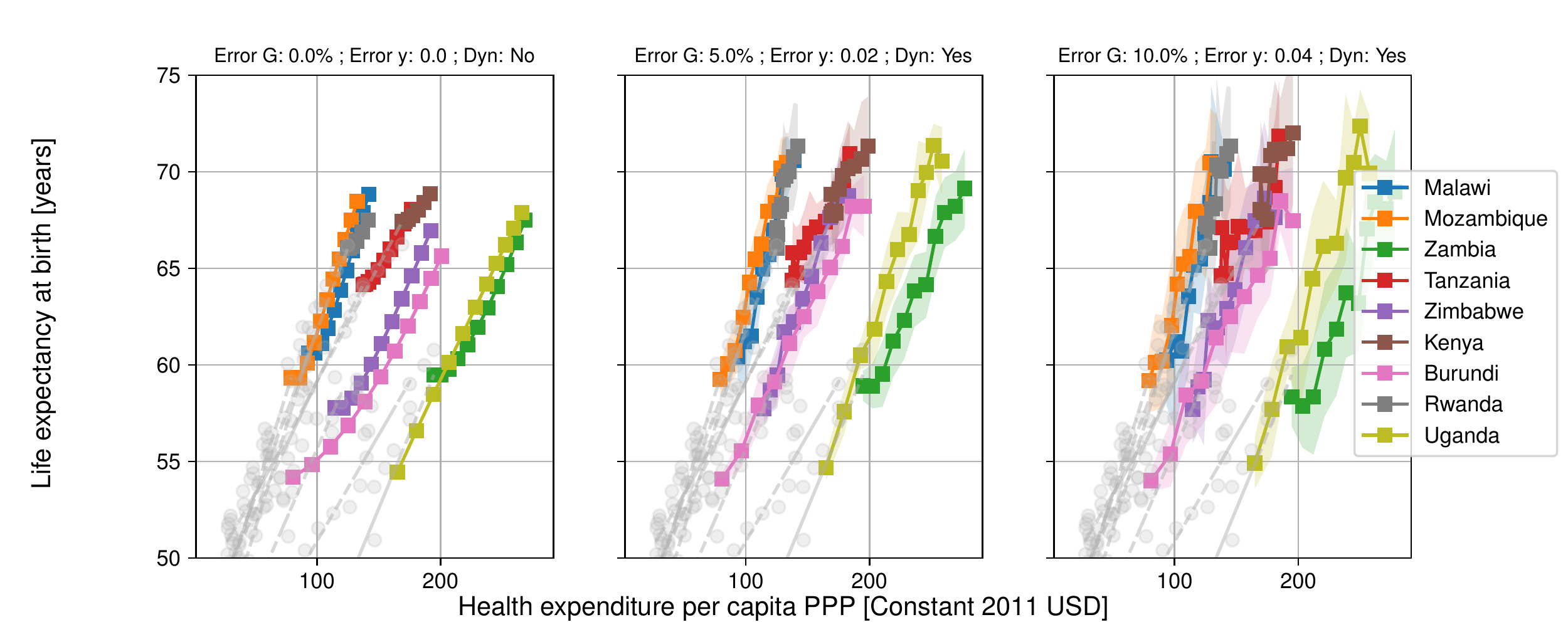}
\caption{Countries' evolution for different scenarios considered in the health subsidies example. When a different dynamics is considered countries can achieve higher or lower life expectancy. }
\label{fig-3}
\end{figure*}

\subsubsection{Results}
In Figure~\ref{fig-2}, we present our first batch of results. For the nine sub-figures, the x-axis represents the time periods (and we generate funding strategies at each period, therefore 10 times), while the y-axis represents the equitability violation metric, which is our cost function. We consider two cases: the lower one represents a static $G_i$; the upper one represents a dynamic $G_i$ as expressed in scenario~(S4). Note that the two sets of curves are in the same scale range, but for reading clarity we moved the upper curves above by using a scaling factor $10$.
In each sub-figure, we represent the solutions we obtain with the three presented policies (SOL, DCL, and DCL+). Each curve is indicated by the mean and the standard deviation over $10$ realizations. 

As we move to the right, we increase the noise error on the feedback $y_i$, from $0.02$ to $0.04$. As we move to the bottom, we increase the error on the erroneously estimated $G_i$ from $5\%$ to $10\%$. Therefore, the top leftmost sub-figure represents the nominal case with and without a dynamic $G_i$, while the bottom rightmost sub-figure is the one with most error sources. 

As we can see, as long as we are in the nominal case with no or small errors, then the SOL policy works quite well, and slightly {outperforms} the other policies (since these perform only one gradient step per time period). However, as soon as the error in the dynamic is higher (as realistic to assume), then both DCL and DCL+ policies are performing better, as one expects. 

Note that DCL and DCL+ are performing very similarly in this example, showing that re-learning $G_i$ at every time step may not be critical in this one-dimensional example, and the feedback is all that one needs.

In Figure~\ref{fig-3}, we plot the trajectories of the different countries when assigned their respective allocations in three different settings. We can see how the dynamics follow closely the previous evolution, and life expectancy increase is slower for countries who are doing already well (e.g., Kenya), and faster for countries that start with a disadvantage (e.g., Uganda). While this is reasonable with our choice of policy (long-term equitability), one can wonder if this is ``fair'' (i.e., \emph{is it ``fair'' to artificially limit better performing countries and invest primarily in worse-off ones?}), and we explore this next.

\subsubsection{Changing cost}

We consider now cases for which we can select a non-zero $\rho$ value in the cost function~\eqref{cost-ex1}. This represents wanting to trade off equal allocation with long-term equitability. For simplicity, we set ourselves in the nominal case and we consider only the SOL policy\footnote{Other policies can be considered, as well as error, but they would not change the main qualitative result.}.

\begin{figure}
\centering
\includegraphics[width = 0.4\textwidth]{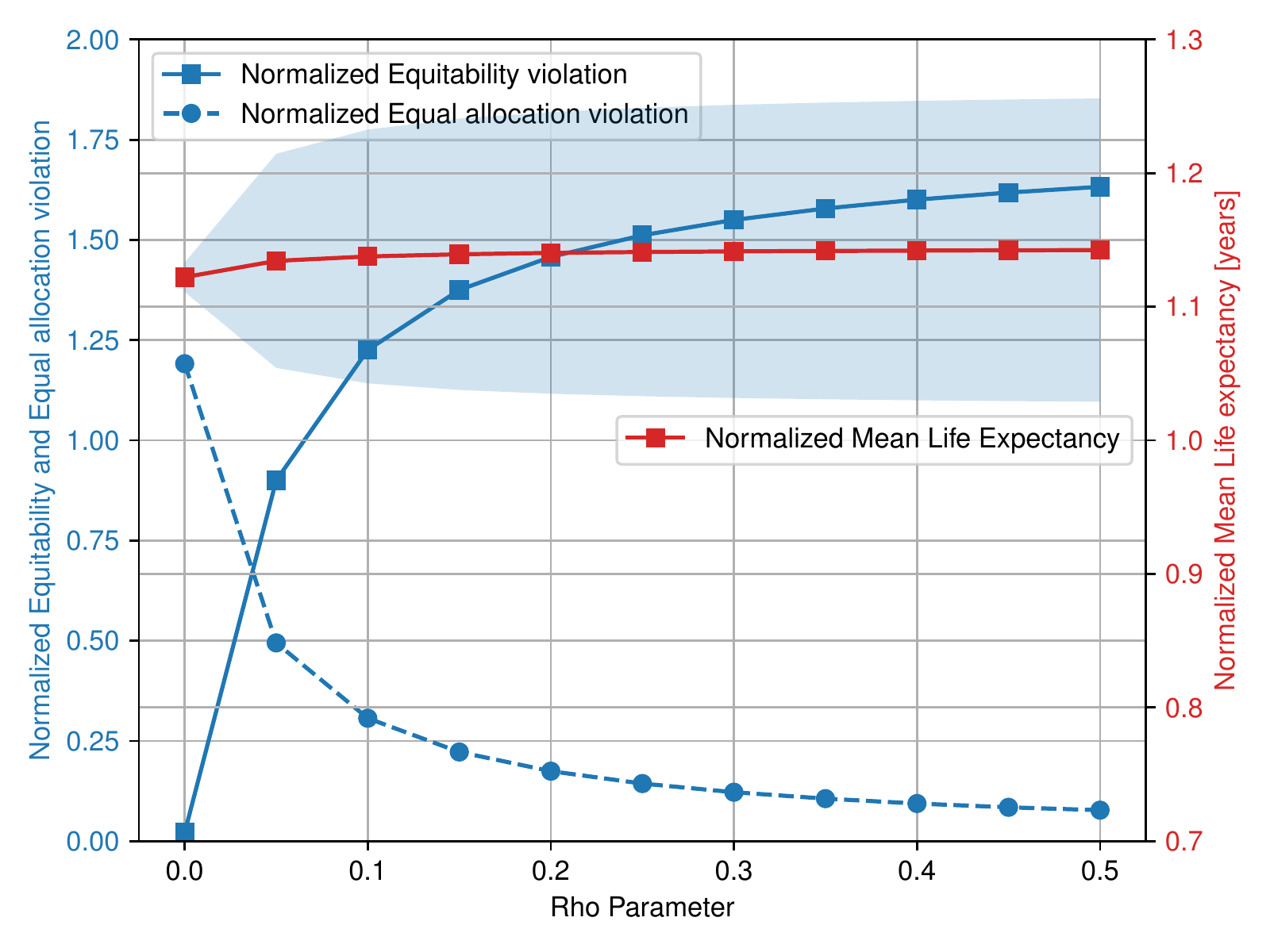}
\caption{Interplay between long-term equitability and equal allocation, in the health subsidies example, for the nominal scenario and SOL policy. Equitability is fragile to compromise (i.e., $\varrho$ values different from $0$), and while the average life expectancy is mainly the same for all values of $\varrho$, the inequalities among countries increase if $\varrho$ increases. This suggests that aiming at imposing ``equality'' by equal allocation can drive inequalities in social outcomes, while not affecting average social welfare. }
\label{fig.3}
\end{figure}

In Figure~\ref{fig.3}, while varying $\varrho$ from $0$ to $0.5$, we plot in blue the equitability violation after ten years normalized with the one at year $0$, and the allocation equality violation after ten years also normalized with the one at year $0$. For these two curves, values below $1$ imply that we are doing better than from when we started from, while values above $1$ signify that we are increasing violations of the metrics. 

We also plot in red the average of the life expectancy after ten years normalized with the one at year $0$, along with its standard deviation. 

As we see, equitability is very sensitive, i.e., \emph{fragile}, to $\varrho$ values different from $0$, and even a small one would increase equitability violation. What is also very remarkable is that despite the average social outcome (i.e., average life expectancy) being mainly the same for all values of $\varrho$, the inequalities among countries increase if $\varrho$ increases, further \emph{segregating} them. This also suggests that aiming at imposing ``equality'' by equal allocation can drive inequalities in social outcomes, while not affecting average social welfare. 

\begin{remark}
This is quite interesting for funding allocation, but also in vaccine allocations over multiple periods, when one wants to make sure that the vaccination uptake is equivalent among neighboring countries (to allow for safe travel), and the global average uptake is as high as possible\footnote{In this case, equitability represents equal uptake among neighboring countries.}. Figure~\ref{fig.3} seems to indicate that the proportional-to-population strategy may not be as effective to achieve \underline{equitability} as more targeted policies based on actual uptakes, and the countries' capacity to turn vaccine vials into vaccinated people. 
\end{remark}

\subsection{Health and Education Subsidies in Malawi}

\subsubsection{Setting}

We look now at a more complex example: the allocation of health and education fundings to a number of local communities in Malawi\footnote{We focus on Malawi, since the country is very active when it comes to education subsidies and education data is available: in particular, data on expenditure corresponds found on World Bank Education Statistics~\cite{dataset-edu-1}, while data on years of schooling found on Barro Lee Education dataset~\cite{dataset-edu-2}.}. The data comes from~\cite{dataset-edu-1, dataset-edu-2}, and we use them to create a realistic, yet fictitious, nominal dynamical system for how the fundings generate welfare and ultimately increases both life expectancy and years of schooling. The nominal system is as follows:
\begin{align*}
  \left[\begin{array}{c} x_{1,k+1} \\ x_{2,k+1} \end{array}\right] 
  & = 
  \left[\begin{array}{cc} 0.5 & 0\\ 0 & 0.3 \end{array}\right]\left[\begin{array}{c} x_{1,k} \\ x_{2,k} \end{array}\right] + \left[\begin{array}{cc} 1.0 & 0\\ 0.01 & 1.0 \end{array}\right] \left[\begin{array}{c} u_{1,k} \\ u_{2,k} \end{array}\right],
\\
\left[\begin{array}{c} y_{1,k} \\ y_{2,k} \end{array}\right] & = \left[\begin{array}{cc}1.0 & 0.03\\ 0.005 & 1.0 \end{array}\right]\left[\begin{array}{c} x_{1,k} \\ x_{2,k} \end{array}\right],
\end{align*} 
where, as in Example~\ref{example}, we set $u_{1,k}$ to be the health funding, while $u_{2,k}$ is the education funding. The output $y_{1,k}$ represents the life expectancy, while $y_{2,k}$ the years of schooling. As we see, the funding in health also generates welfare in education. The education dynamic $0.3$ has less inertia than the health, which is also empirically observed. This system is compatible with the observed data and it will be used to evaluate our policies (which, we remind, do not need access to the system to work, just input-output data). 

We generate $25$ local communities by adding noise to the coefficients of the above nominal system\footnote{In particular, we add zero-mean Gaussian noise terms with different standard deviation $S$: $S_{A_{11}} = 0.02$, $S_{A_{22}} = 0.05$, $S_{B_{21}} = 0.0025$, $S_{C_{12}} = 1.5$e${-5}$, $S_{C_{21}} = 2.5$e$-4$.}, and generate a random network of interconnection to model graph $\mathcal{G}$. We run our policies on a scenario in which we want to decide the allocation of the fundings for the next $10$ funding periods (years) in order to optimize long-term equitability. We assume we have no feedback noise, but we insert a $5\%$ random noise level on the estimation of $G_i$, and consider $25$ different realizations. As in the first example, we consider increasing the funding at each period (so to have $50\%$ for health, and $25\%$ for education in ten years). We consider cost function~\eqref{cost-ex1} with $\varrho = 0$, to focus on long-term equitability only.  

\subsubsection{Results}

\begin{figure}
\centering
\hspace*{-0.75cm}
\includegraphics[width = 0.575\textwidth]{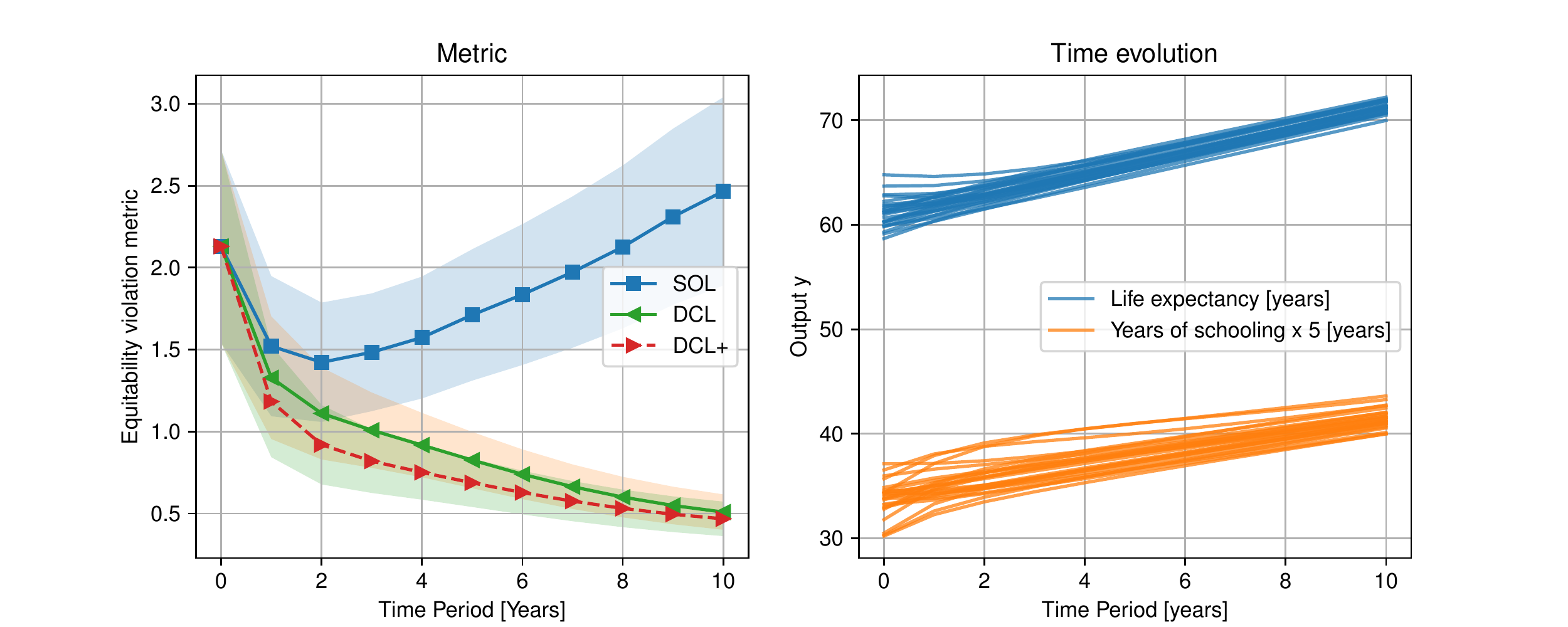}
\caption{Metric and evolution for the example of subsidies for health and education in $25$ local communities. Here we see how the DCL+ policy behaves the best in reducing equitability violation. }
\label{fig.4}
\end{figure}

In Figure~\ref{fig.4}-left, we can see how the policies contribute at diminishing the violation of long-term equitability. As we can see, since we have an estimation error on the estimation of $G_i$, we observe that policy SOL behaves the worst, while integrating feedback is the best solution. As in the previous example, we plot here the mean and we shade the standard deviation along the mean for the different realizations. 

We notice how policy DCL+ performs the best in this two-dimensional scenario, advocating for its need, in general, for more complex dynamical systems with correlations among the different dimensions.  

Figure~\ref{fig.4}-right depicts (for one realization) the evolution of both the life expectancy and the years of schooling for the different communities. As we see, the respective values are driven to coalesce.  

\subsubsection{Changing cost}

\begin{figure}
\centering
\includegraphics[width = 0.45\textwidth]{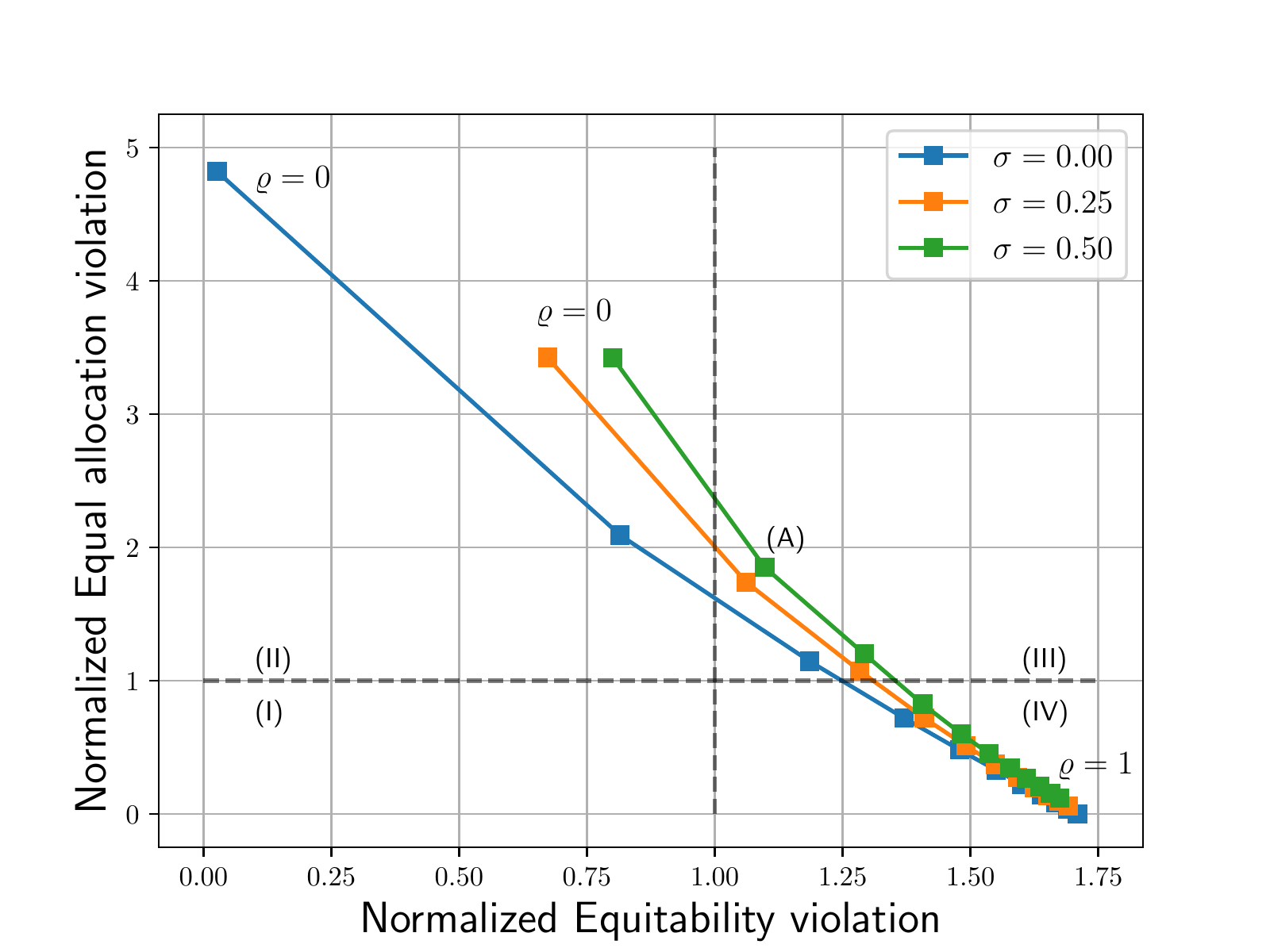}
\caption{Pareto frontiers for the example of subsidies for health and education in $25$ local communities. This figure represents the trade-offs and compromises inherent in any democratic country: balancing local preferences (big $\sigma$), equitability, and allocation equality (big $\rho$) in different ways lead to different social outcomes. The main takeaway from this figure is that a naive compromise could lead to point (A), where one degrades both equal allocation and equitability. Also equitability is fragile to compromise (as we have comparably less points in quadrant (II) ).
}
\label{fig.6}
\end{figure}

In Figure~\ref{fig.6}, we expand the setting. For simplicity, we consider no estimation error, so SOL behaves well, and we do not have issues related to different realizations, while maintaining the same qualitative result. Then we consider a more complex cost function, as
\begin{equation}\label{cost-ex2}
f_{\textrm{ex.2}} =  \!\!\sum_{i,j \in [1,N] }\varrho\Big \| \, \bar{u}_i - \bar{u}_j \,\Big\|^2 + \!\!\sum_{i\in [1,N] }\Big[(1-\varrho)\Big\| \, \bar{y}_i - \frac{1}{N_i}  \sum_{j \in \mathcal{N}_i}\bar{y}_j \,\Big\|^2 + \sigma \Delta_i\Big],
\end{equation} 
trading-off equality of allocation, long-term equitability, and communities' preferences. In particular, we set: 
\begin{equation}
\Delta_i = \omega_i^u \Big\| \, \bar{u}_i - \frac{1}{N_i}  \sum_{j \in \mathcal{N}_i} \bar{u}_j \,\Big\|^2 + \omega_i^y \Big\| \, \yik - \frac{1}{N_i}  \sum_{j \in \mathcal{N}_i}\yjk \,\Big\|^2,
\end{equation}
where $\omega_i^u, \omega_i^y \in [0,1]$ are non-negative weights that capture each community preference to equality of allocation or long-term equitability within their neighborhood. We set $\omega_i^u = 0, \omega_i^y = 1$ for the first $13$ communities and the other way around for the remaining $12$. Note that $\omega_i^u, \omega_i^y$ could be learned on-line as explained in e.g.,~\cite{Simonetto2021, Notarnicola2022}, but we do not do it here for the sake of simplicity. 

In Figure~\ref{fig.6}, we capture the results in terms of Pareto frontier. In particular for different choices of $\sigma = \{0, 0.25, 0.5\}$, we vary $\varrho$ from $0$ (all long-term equitability) to $1$ (all equal allocations), and we plot the normalized (to the initial time period $k=0$) end-value of equitability violation (i.e, value at $k=10$/ value at $k=0$) as well as normalized end-value of equal allocation violation. 

We observe the following. Consider the case $\sigma = 0$, so no community personalization is present. When $\varrho$ varies smoothly from $0$ to $1$, then the equitability violation gets worse, and equal allocation violation gets better, as one expects. We find back the fact that long-term equitability is fragile, since when $\varrho$ is slightly $>0$ (here $\varrho = 0.2$) then equitability violation gets $>1$, and therefore worse than the beginning. Equal allocation is less fragile, seemingly having an accumulation of points for large values of $\varrho$. Observe also the point for which both equitability and equal allocation are both $>1$, signifying that both metrics get worse than the beginning. This possibility should not surprise, since the initial condition is not at equilibrium, it is possible for things to get worse. Which happens if one is aiming at a difficult balance between two conflicting metrics, without a decisive direction.   

Consider now the case $\sigma > 0$, having more than half of the communities wanting long-term equitability, and the rest equal allocation. We see again the observed trend for $\sigma=0$, with the important difference that long-term equitability violation is larger. We find again the possibility for both equitability and equal allocation to get worse than when we started (both values $>1$), which should not surprise, but it must be avoided in practice. 

Focus on point (A) in the graph. There the governmental funding agency decides to weigh in the local communities' feedback, while balancing equitability with equal allocation: from many perspectives, its policy is a justified compromise. However, things get worse for both metrics. Even if for some of the communities a local metric may have improved, as a whole the country is worse off, suggesting that (naive) compromise may not be always a good strategy\footnote{This does not mean that compromise should be avoided at all costs. Our results suggest that compromise should be planned carefully, taking into account the community dynamic evolutions, as we do here. }. 

Observe again Figure~\ref{fig.6}, and in particular how few points are in the quadrant (II) that improves equitability: equitability can be achieved, but it is fragile. It is very sensitive to compromise and trade-offs. 

We close looking at quadrant (I): the region where everything improves. While one would like to drive the country there, this is independent from subsidies. As we have argued in Section~\ref{sec:design}, the only possibility to have an equilibrium in this area is that $G_i \approx G_j$, so that equal allocation implies equitability and vice versa. This further suggests that subsidy design alone is not sufficient to drive countries and communities to increase welfare, but one also needs systemic changes to how the welfare is generated and transformed (i.e., one needs to change the underlying dynamical systems, and therefore $G_i$).

\section{Conclusions}

We have studied long-term equitability for allocating limited resources in a dynamic setting, whereby local communities can evolve based on the allocations that a funding agency provides, and they give feedback on their social outcome, and possibly preferences. We have proposed policies to drive the system to the desired long-term equitability and we have empirically shown how this equitability is fragile. In particular, it can be quickly lost if allocations are divided to take also into account other fairness notions, such as funding parity. 

As such, the main takeaways of this paper are as follows: 

$\bullet$ Whenever one considers dynamic effects and looks at the interplay between decisions and how who is affected by them evolves because of them, then incorporating feedback is critical in designing policies that are robust to modeling errors. These online optimization with feedback policies can then achieve long-term equitability, if the cost function is properly tuned. 

$\bullet$ Long-term equitability is very sensitive to compromise. If one decides to balance local community wishes with global equitability and other forms of parity (e.g., funding parity), then any gains in equitability can be quickly lost, and naive compromises can even lead to the worsening of most (or all) of the social outcomes that were considered in the compromise. 

These two conclusions should be central whenever allocating limited resources in a dynamic setting, ranging from designing policies for health subsidies (in low and high income countries alike) to designing subsidies for the transition to a net-zero carbon world.


\bibliographystyle{acm}
\balance
\bibliography{references}

\end{document}